%% file: duminil-bio2022.tex
\documentclass[12pt]{amsart}

\RequirePackage[colorlinks,citecolor=blue,urlcolor=blue]{hyperref}
\usepackage{verbatim}
\usepackage{graphicx}

\newcommand\Tc{T_{\text{\rm c}}}
\newcommand\resp{respectively}
\newcommand\ZZ{{\mathbb Z}}
\renewcommand\o{\text{\rm o}}

\newcommand\DC{Duminil-Copin}

\begin{document}

\title{The work of Hugo Duminil-Copin
}
\author{Geoffrey R. Grimmett}
\address{Statistical Laboratory, Centre for
Mathematical Sciences, Cambridge University, Wilberforce Road,
Cambridge CB3 0WB, UK} 
\email{grg@statslab.cam.ac.uk}
\urladdr{\url{http://www.statslab.cam.ac.uk/~grg/}}
\date{\today}

\begin{abstract}
This article is an account of the scientific work of Hugo Duminil-Copin at the time of his award
in 2022 of the Fields Medal \lq\lq for solving longstanding problems in the probabilistic theory of 
phase transitions in statistical physics, especially in dimensions three and four''.
\end{abstract}

\maketitle

\centerline{\includegraphics[width=0.75\textwidth]{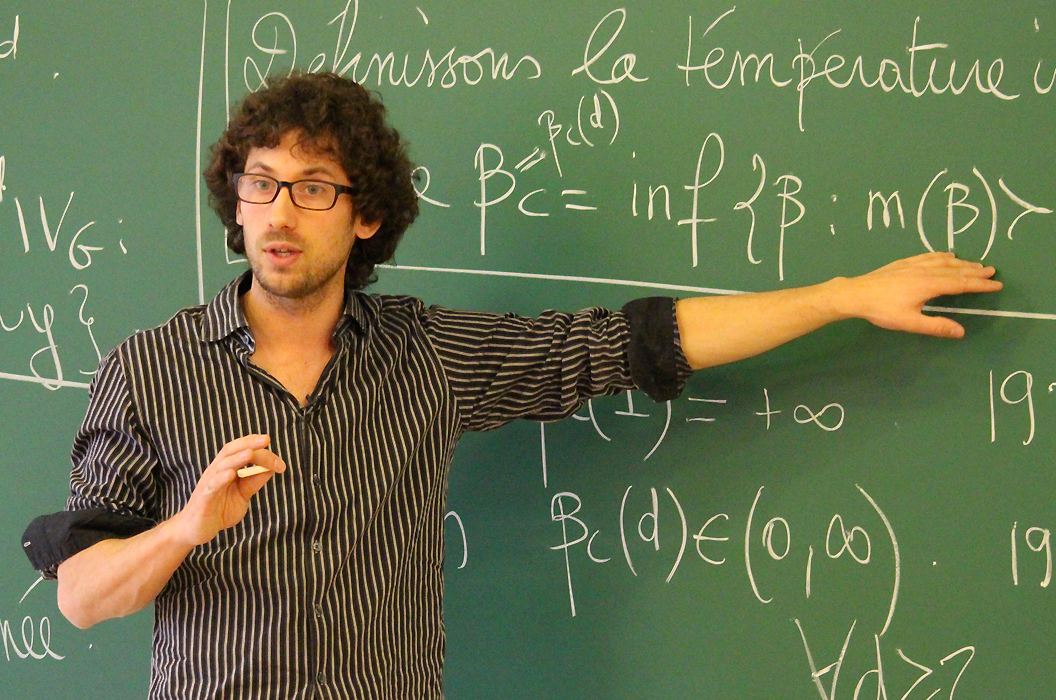}}

\section{Overview}

{\sc Magnetism.}
Ice melts, water boils --- these are examples of \lq phase transitions'.
A further phase transition 
was studied experimentally by Pierre Curie for his 1895 Doctoral Thesis.
A magnetized bar of iron retains its magnetization if and only if
the temperature $T$  is sufficiently low, or more specifically $T\le \Tc$ for a
 critical temperature $\Tc$ now termed 
 the \lq Curie point'. (In the case of iron, we have $\Tc\approx 770^\circ$C.)

\smallskip
{\sc Spatial disorder.} Gas particles are filtered through a disordered material.
On a microscopic scale, the medium comprises fissures some but not all of
which are wide enough to allow the gas to pass.  
The practical question is whether or not the gas passes through the filter.

\begin{figure}
\centerline{\includegraphics[width=0.65\textwidth]{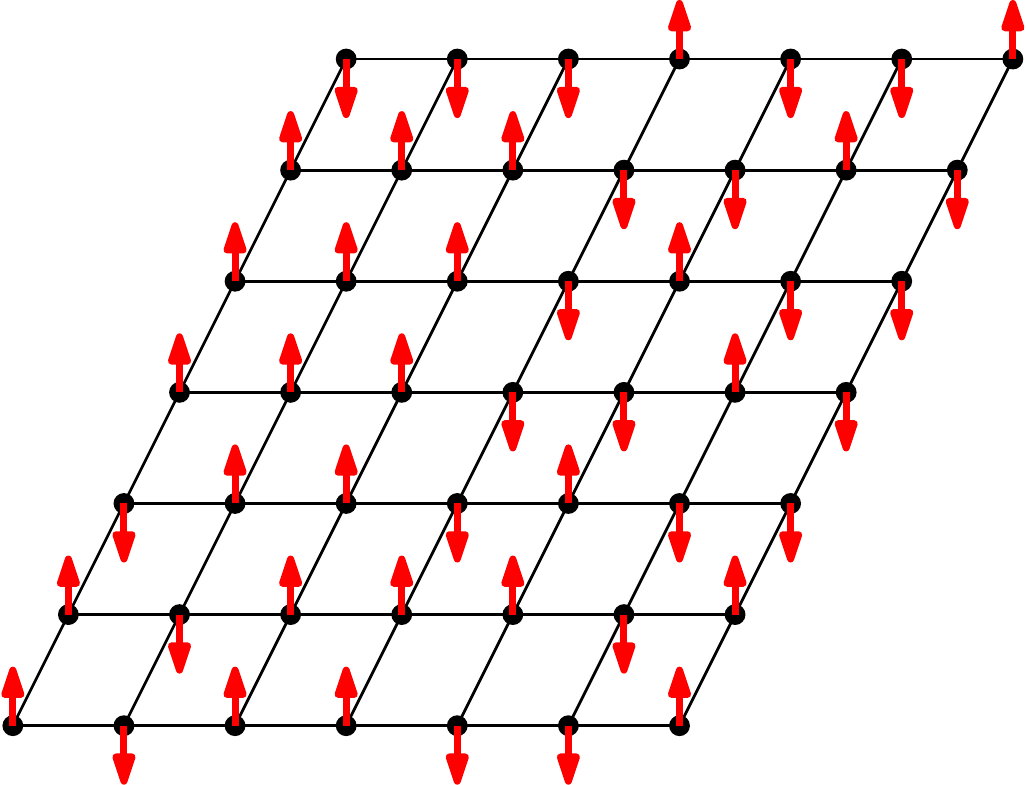}}
\caption{A representation of the Ising model on the square lattice. Each vertex represents
a particle that has a charge modelled as either a north pole or a south pole. The probability of a configuration
depends on the number of nearest-neighbour pairs with like poles.}
\label{fig:ising}
\end{figure}

\begin{figure}
\centerline{\includegraphics[width=0.65\textwidth]{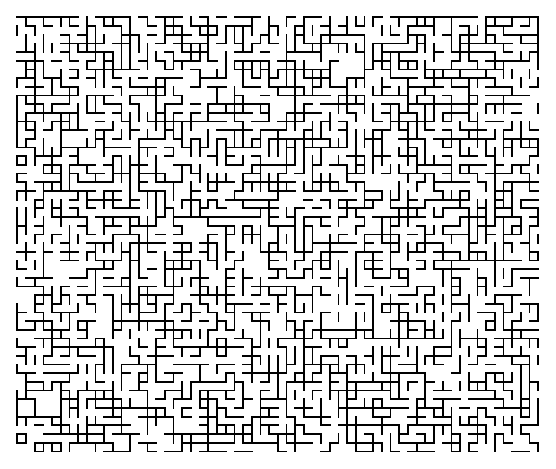}}
\caption{A simulation of bond percolation on $\ZZ^2$ in which each edge is designated \lq open' wth
probability with $p=0.51$. This process is barely supercritical.}
\label{fig:perc}
\end{figure}

\smallskip
Mathematicians and theoretical physicists study such systems by building 
and analysing conceptual models. The fundamental models for magnetism and spatial
disorder are, \resp, the \lq Ising model' (of 1920) and the \lq percolation model' (of 1957).
See Figures \ref{fig:ising} and \ref{fig:perc}.
These two models coexist within a large  (and growing) family of models for disorder
in physical systems. Together with their cousins, they have 
attracted an enormous amount of attention from mathematicians and physicists. 
Each model possesses a phase transition, and
the principal challenges lie in understanding the nature of this transition.
The mathematics of such systems is ramified, highly technical, and  very complex.
Many significant problems have been overcome, and many remain. For example,
the Ising model in two dimensions is now largely understood by mathematicians, 
whereas physicists are well ahead of mathematicians in three dimensions.

In recent years, no-one has advanced the mathematical theory of these and related
systems more than Hugo \DC. 
Hugo is an extraordinary mathematician who has achieved enormous progress. 
His many 
contributions to the theory of critical phenomena have resolved 
longstanding problems of deep significance, 
and have opened new directions of research. 

From his 2011 PhD thesis onwards, \DC\ has had unique impact on the rigorous theory. 
In a  large number of publications characterized by
both depth and breadth, he has solved a wide range of important problems, and has
introduced novelty and clarity to the methodology. At the same time, he has displayed outstanding
qualities of leadership through his mentoring of early-career researchers. 

\DC's early work was based largely around the idea that critical systems in two dimensions 
are invariant under so-called \lq conformal mappings'; these are the mappings that
are composed of local dilations and rotations. This topic has 
seen an explosion of progress over the last twenty years, 
starting with Oded Schramm's introduction of a type of random curve termed SLE.
While still in his twenties, \DC\ made several fundamental contributions,
including rigorous derivations of long-conjectured critical probabilities and connective constants.

The hypothesis of \lq universality' asserts in such contexts that the nature of a phase 
transition depends only on the type of model and the number of dimensions; thus, for example, the phase transitions of percolation on the square and triangular lattices
are expected to be of common type.   
\DC\ expanded his research interests to the study of universality in statistical physics, where he obtained a number of outstanding results, including the scaling relations for 
certain random-cluster models that incorporate both Ising and percolation models, 
together with fractal properties of continuous phase transitions. 

Moving beyond two dimensions, \DC\  has devised new methods
to resolve a number of long-open and notorious classical problems in critical 
phenomena, including: (i)
the    continuity of the phase transition of the three-dimensional Ising model,
(ii) the `triviality' of the Ising model scaling limits in four dimensions, and 
(iii) the sharpness of the phase transition for a range of stochastic models in arbitrary dimension. 
He has also made many significant and universal advances in the theory of off-critical phases
of a number of statistical physics models of importance.

In the following sections, we describe in more detail  a selection of \DC's 
profound and diverse contributions.
For summaries of some of his results, the reader may consult his expository work
\cite{expo-sbm,cdm,ecm, icm,expo-jpn,expo-ising}.

Most of Hugo's papers have been written jointly with others.  
He is invariably enthusiastic to share his ideas openly, and he
collaborates freely and intensively with 
an unusually wide spread of colleagues, both younger and older. 
This practice has greatly augmented and amplified his beneficial impact on science. 

\section{Percolation, Ising, Potts models}

The random-cluster model was introduced by Fortuin and Kasteleyn around 1970 as a
unification of electrical networks, percolation, and the Ising/Potts models. It has become
a focus of unification in the theories of these topics.

It was long conjectured that the critical point of the
random-cluster model on the square lattice with cluster-weight $q$
is the self-dual point 
$$
p_c(q) = \frac{\sqrt q}{1+\sqrt q}.
$$
The $q=1$ case is a famous result of
Harris and Kesten (1980) for percolation, and the $q=2$ case amounts
to an Onsager calculation (1944) for the Ising model.
Some progress had been made by others for certain ranges of $q$, 
but the general result remained open until proved
in 2012 by Duminil \cite{pc-beffara} (with Beffara). 
Their solution resolved, in particular, the (longstanding) conjectured value of the critical point  of the 
$q$-state Potts model on the square lattice.

In a  subsequent work, \cite{sharp-raoufi-tassion}, to which we return later, he 
introduced a new method for proving exponential decay (and more),
and applied it to a spectrum of processes including all random-cluster models and Potts models. This has proved a major step towards a systematic theory of disordered systems.

Duminil-Copin has been instrumental in extending many of the classical percolation techniques from  
Bernoulli percolation (independent) models to a much wider class of dependent models. 
With Hongler and Nolin \cite{rsw-hongler-nolin}, and later Chelkak and Hongler 
\cite{rsw-chelkak-hongler}, 
he proved  RSW-type \emph{a priori} bounds for probabilities of crossings
in the FK Ising model. This was an early example of an RSW-type theory for percolation models with
dependence, and it paved the way to many later advances by Duminil-Copin and other authors.

Much attention has been given in the physics literature to the question of (dis)continuity
in the phase transition of the two-dimensional Potts model. This was resolved by \DC\ 
in \cite{discontinuity} (with Gagnebin, Harel, Manolescu,
Tassion) with the proof of discontinuity when $q>4$, and in \cite{continuity}
(with Sidoravicius and Tassion) with
continuity when $1\le q \le 4$.
  
 \begin{figure}
 \centerline{\includegraphics[width=0.6\textwidth]{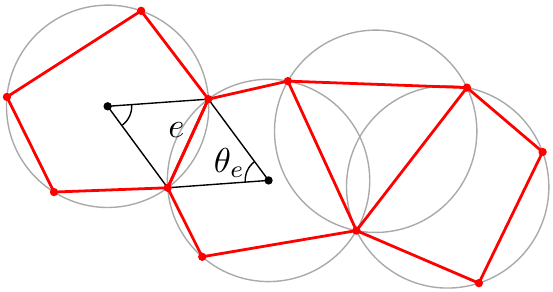}}
  \caption{A plane graph is called \emph{isoradial} if
  its faces have circumcircles of constant radius.}
  \label{fig:iso}
 \end{figure}
  
In two outstanding recent papers, with Li and Manolescu \cite{iso-li-manolescu} and 
Manolescu  \cite{scaling-manolescu}, \DC\  has proved certain scaling properties of
random-cluster models. The first paper concerns isoradial graphs, as illustrated in Figure \ref{fig:iso}.
Isoradiality enables the identification of the canonical critical points predicted by the star--triangle transformation
(or Yang--Baxter relation); by judicious use of the star--triangle transformation,
certain scaling properties may be propagated across the family of critical isoradial processes.
The second paper establishes scaling relations for the critical model on the square lattice.  
Such results were known before only in the percolation case $q=1$
in famous work of Kesten.
This pair of papers solve several big open problems, and constitute 
key progress towards proving universality 
for random-cluster models on isoradial graphs.

One of \DC's recent results in this area is the striking proof \cite{rotational}
(joint with  Kozlowski, Krachun, Manolescu, Oulamara)
 that the model on the square lattice with parameter $q\in [1,4]$ exhibits \emph{rotational invariance} at large scales. 
The proof is a beautiful and technically demanding combination of several methods, including
modification of the lattice by star--triangle transformations and Bethe ansatz calculations.
This impressive development
opens the way to promising approaches for proving full conformal invariance.

\section{Self-avoiding walks}

\begin{figure}
 \centerline{\includegraphics[width=0.4\textwidth]{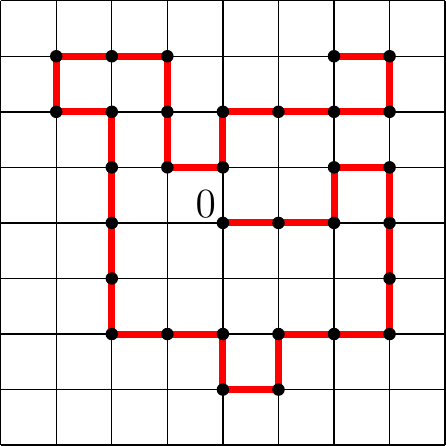}}
  \caption{A self-avoiding walk of length 31 on the square lattice, starting at the origin.}
  \label{fig:saw}
 \end{figure}

The theory of self-avoiding walks  (SAWs) offers perhaps some of
the most challenging open problems in probability and combinatorial theory. 
(See Figure \ref{fig:saw}.) 
\DC\ provided a startling solution to a notorious conjecture that emerged
in 1982 work of Nienhuis in conformal field theory. He proved in \cite{saw-smirnov} 
(joint with Smirnov)  that
the number of $n$-step self-avoiding walks on the hexagonal lattice grows
in the manner of  $\mu^{n+\o(n)}$ with
$\mu = \sqrt {2 + \sqrt {2}}$.
The  elegant proof uses a parafermionic observable (similar to the holomorphic
observable used by Chelkak and Smirnov in the context of the Ising model),
together with some neat counting/convergence arguments. 
This major result has opened the way to proving convergence to SLE of a random SAW
in two dimensions.

He has obtained a number of other important results
about SAWs. Firstly,  he showed in \cite{saw-hammond} (with Hammond)
that SAW in dimensions $2$ and more are 
sub-ballistic; this may be viewed as the
first non-trivial result in this direction for the `difficult'
dimensions 2 and 3. 
Secondly, in the other direction his papers \cite{saw-glazman-hammond-manolescu}
(with Glazman, Hammond and Manolescu) concerning delocalization, 
and \cite{MR3189073} (with Kozma and Yadin) concerning the space-filling property of  
supercritical SAW, are notable.

\section{Ising and $\phi^4_4$ lattice models}

\DC\  has addressed several classical problems of critical 
phenomena which first achieved prominence in the 1980s.  
 He began this project by identifying a list of fundamental challenges, being typically 
 fundamental  questions which were still viewed as being out of reach.  
 He then achieved these goals
 in sequence.  In some cases this was done via  bold enhancements  of 
 known methods, 
 and  in others through the introduction of  new ideas.   

He showed in \cite{ising-aizenman-sidoravicius,ecm} (with Aizenman and Sidoravicius)
how to develop the random-current representation of the Ising model to resolve the 
problem of proving the continuity of the Ising phase transition in three dimensions.
This is perhaps the first such result for any classical system in three dimensions, 
and it completes
the picture for the Ising model in general dimensions.
The proof techniques have since been utilized by several colleagues.

In the paper \cite{phi-aizenman} (with Aizenman), Duminil-Copin
showed the \lq marginal triviality\rq\ of the scaling limits of Ising and $\phi^4_4$ lattice models in four dimensions.
Namely, the scaling  limits of (near-)critical Ising and $\phi^4_4$ models 
are trivial (or Gaussian) as Euclidean field theories. 
Such triviality was derived in dimensions $d \ge5$ by Aizenman and Fr\"ohlich 
(independently) in 1981/82.
The recent extension to the (critical) dimension $d=4$  was obtained 
via a singular correction to previously deployed bounds. 
The last 
 was derived through a multiscale analysis of the random-current intersection probability 
 at the threshold dimension 4. 

\section{Sharpness of phase transitions}

A major thrust of \DC's work has been to understand the sub- and supercritical 
phases of models of statistical physics. There are classical methods that
allow study of small and large parameter-values, but it is very much harder to
study the phases right to the critical point. \DC\ has many influential results in this direction,
beginning with his proof of the critical point of the two-dimensional random-cluster model,
\cite{pc-beffara}.

In \cite{sharp-raoufi-tassion}, \DC\ (with Raoufi and Tassion) developed a beautiful method to prove exponential decay in the subcritical phase,
therein making use of an inequality on decision trees of 
O'Donnell et al. The authors explained, in particular,
how to use this for random-cluster models
in all dimensions. This last result, 
which was previously known only when either $q=1$ or $q$ is sufficiently large, 
answers a problem of about 30 years standing. 
The method has many other novel
applications, including the best proof of exponential decay
for Voronoi percolation,
\cite{MR3916112}. For a two-dimensional system with a property of duality, it 
leads to
the coincidence of the self-dual point and the critical point.
This work has changed our vision of 
the property of so-called \lq sharp threshold\rq.  

One important and topical application of the above method
is given in \cite{gff-goswami-etal} (with Goswami, Rodriguez and Severo), where \DC\ has
studied phase transition for the $h$-level set of the Gaussian free field in dimension $\geq 3$,
and has proved the sharpness of the critical value of $h$.

\section{Further work}

Hugo has a number of outstanding contributions that do not fit easily into 
the above classification, of which a few examples are mentioned here.

He established in  \cite{bootstrap}
the sharp threshold for bootstrap percolation in all dimensions 
 (with Balogh,  Bollob\'as and Morris).

In \cite{percolation-goswami-etal} (with   Goswami,  Raoufi,  Severo and  Yadin),
he established the Benjamini--Schramm conjecture that Bernoulli percolation on a
Cayley graph with superlinear growth has $p_c<1$.
This was achieved by
exploiting  a new connection between percolation and the Gaussian free field, which
relates the connectivity properties of percolation clusters to the geometry of the graph.
    
Paper \cite{dim-benjamini} (with Benjamini, Kozma and Yadin) is an interesting work on the
dimension of spaces of harmonic functions on certain random graphs. This issue is
connected with a number of other topics including the theory of  random walks in random environments.

\section{Summary}

Hugo \DC\ has found new ways of looking at old problems. He has  obtained solutions
to classical problems of great visibility,
while simultaneously introducing new methods and intuitions into the field.
His area of mathematical science has been largely reshaped by his achievements.

Hugo displays the strongest personal attributes in addition to his scientific talent. 
The fact that almost all his papers
are in collaboration reflects his desire to communicate with and to excite others,
and he invariably succeeds.

\section*{Ackowledgements}
The writer acknowledges  the significant contributions of a number of colleagues to this scientific
biography, namely Michael Aizenman, Erwin Bolthausen, Stanislav Smirnov, and Ofer Zeitouni. 
The photograph is included with the kind permission of IHES.

\input{duminil-bio2022.bbl}

   \end{document}

%% file: duminil-bio2022.bbl
\providecommand{\bysame}{\leavevmode\hbox to3em{\hrulefill}\thinspace}
\providecommand{\MR}{\relax\ifhmode\unskip\space\fi MR }
\providecommand{\MRhref}[2]{%
  \href{http://www.ams.org/mathscinet-getitem?mr=#1}{#2}
}
\providecommand{\href}[2]{#2}